\documentclass{article}

\usepackage{amsmath}
\usepackage{amsthm}
\usepackage[psamsfonts]{amsfonts}
\usepackage{amssymb}
\usepackage[all]{xy}

\DeclareMathSymbol{\rightrightarrows}  {\mathrel}{AMSa}{"13}

\def\Aut{\operatorname{Aut}}

\def\Ob{\operatorname{Ob}}
\def\Mor{\operatorname{Mor}}

\def\St{\operatorname{St}}

\def\Ho{\operatorname{Ho}}

\def\Iso{\operatorname{Iso}}

\catcode`\@=11\def\varholim@#1#2{\mathop{\vtop{\ialign{##\crcr
 \hfil$#1\m@th\operator@font holim$\hfil\crcr
 \noalign{\nointerlineskip\kern\ex@}#2#1\crcr
 \noalign{\nointerlineskip\kern-\ex@}\crcr}}}}
\def\hocolim{\mathpalette\varholim@\rightarrowfill@} 
\def\hoinvlim{\mathpalette\varholim@\leftarrowfill@}
\catcode`\@=\active

\newtheorem{theorem}{Theorem}
\newtheorem{lemma}[theorem]{Lemma}
\newtheorem{corollary}[theorem]{Corollary}

\theoremstyle{definition}

\newtheorem{example}[theorem]{Example}

\newtheorem{remark}[theorem]{Remark}

\begin{document}

\title{Homotopy classification of gerbes}
 
\author{J.F. Jardine}
 
 
\date{March 27, 2006}
 
\maketitle

\begin{abstract}
\noindent
Gerbes are locally connected presheaves of groupoids on a small
Groth\-endieck site $\mathcal{C}$. Gerbes are classified up to local
weak equivalence by path components of a cocycle category taking
values in the diagram $\mathbf{Grp}(\mathcal{C})$ of $2$-groupoids
consisting of all sheaves of groups, their isomorphisms and
homotopies. If $\mathcal{F}$ is a full subpresheaf of
$\mathbf{Grp}(\mathcal{C})$ then the set $[\ast,B\mathcal{F}]$ of
morphisms in the homotopy category of simplicial presheaves classifies
gerbes locally equivalent to objects of $\mathcal{F}$ up to weak
equivalence. If $\St(\pi \mathcal{F})$ is the stack completion of the
fundamental groupoid $\pi\mathcal{F}$ of $\mathcal{F}$, if $L$ is a
global section of $\St(\pi\mathcal{F})$, and if $F_{L}$ is the
homotopy fibre over $L$ of the canonical map $B\mathcal{F} \to
B\St(\pi\mathcal{F})$, then $[\ast,F_{L}]$ is in bijective
correspondence with Giraud's non-abelian cohomology object
$H^{2}(\mathcal{C},L)$ of equivalence classes of gerbes with band $L$.
\end{abstract}

 \section*{Introduction}

Suppose that $\mathcal{M}$ is a closed model category, and that $X$
and $Y$ are objects of $\mathcal{M}$. A cocycle from $X$ to $Y$ is a
picture
\begin{equation*}
X \xleftarrow{f} Z \xrightarrow{g} Y
\end{equation*}
of morphisms in $\mathcal{M}$ such that $f$ is a weak equivalence. A
morphism of cocycles $(f,g) \to (f',g')$ is a commutative diagram
\begin{equation*}
\xymatrix@R=8pt{
& Z \ar[dl]_{f} \ar[dd] \ar[dr]^{g} & \\
X && Y \\
& Z' \ar[ul]^{f'} \ar[ur]_{g'}
}
\end{equation*}
and these cocycles and their morphisms together form the category
$H(X,Y)$ of cocycles from $X$ to $Y$. 
The assignment $(f,g) \mapsto gf^{-1}$ defines a function
\begin{equation*}
\phi: \pi_{0}H(X,Y) \to [X,Y]
\end{equation*}
from the path components of cocycle category $H(X,Y)$ to the set of
morphisms $[X,Y]$ from $X$ to $Y$ in the homotopy category
$\Ho(\mathcal{M})$. Then it is a basic result of \cite{J4} that this
function $\phi$ is a bijection if the model category $\mathcal{M}$ is
right proper and if its class of weak equivalences is closed under
finite products.

The right properness condition is a serious restriction, but right
proper model structures are fairly common in nature, and include the
standard model structures for spaces, simplicial sets, and spectra, as
well as more exotic structures such as simplicial presheaves,
simplicial sheaves and presheaves of spectra on small Grothendieck
sites. 

The cocycle approach to constructing morphisms in the homotopy
category is proving to be very useful, particularly in connection with
simplicial sheaves and presheaves. Applications have so far appeared
in new, short and conceptual arguments for the homotopy classification
of sheaf cohomology theories, both abelian and non-abelian
\cite{J4}. Cocycle categories are involved in the explicit
construction of the stack completion functor which is given in
\cite{J3}. They have been used to show \cite{J2}, in a variety of
settings, that morphisms $[\ast,BI]$ in the homotopy category can be
identified with path components of a suitably defined category of
$I$-torsors for all small category objects $I$.

The present paper uses cocycles in presheaves of $2$-groupoids, here
called $2$-cocycles, to give a homotopy classification of gerbes.

A gerbe is typically defined in the literature \cite[p.129]{Gir} to be
a stack $G$ which is locally path connected. Stacks themselves have no
conceptual mystery: they are fibrant objects in model structures for
sheaves of groupoids \cite{JT} or more generally presheaves of
groupoids \cite{H}, and one now can identify a stack with the homotopy
type that it represents in presheaves of groupoids. The model
structure for presheaves (or sheaves) of groupoids, over any small
Grothendieck site $\mathcal{C}$, is easy to describe: a map $f: G \to
H$ of presheaves of groupoids is a local weak equivalence
(respectively global fibration) if the induced map $BG \to BH$ of
classifying objects is a local weak equivalence (respectively global
fibration) of simplicial presheaves. Local path connectedness is an
invariant of homotopy type in this sense, and we shall take the point
of view that a gerbe is a presheaf of groupoids $G$ such that the
classifying simplicial presheaf $BG$ is locally path connected. The
local path connnectedness condition can be expressed this way: given
objects $x,y \in \Ob (G)(U)$ in a section $G(U)$, there is a covering
family $\phi: V \to U$ such that there is a morphism $\phi^{\ast}(x)
\to \phi^{\ast}(y)$ in $G(V)$ for any $\phi$ in the cover.

Every gerbe $G$ is locally equivalent to any of its sheaves of
automorphism groups. The category $H(\ast,\mathbf{Grp}(\mathcal{C}))$
of $2$-cocycles taking values in the diagram of all sheaves of groups,
their isomorphisms and homotopies, is the vehicle by which we classify
gerbes up to local weak equivalence. Theorem \ref{thm 20} says
that the path components of this cocycle category are in one to one
correspondence with the path components of the category
$\mathbf{Gerbe}(\mathcal{C})$ of gerbes and their local weak
equivalences, or that there is a bijection
\begin{equation*}
\pi_{0}(\mathbf{Gerbe}(\mathcal{C})) \cong 
\pi_{0}H(\ast,\mathbf{Grp}(\mathcal{C})).
\end{equation*}
One has to interpret a statement like this correctly, because the
categories involved are not small. The path component functor
$\pi_{0}$ means the class of equivalence classes of objects, where
two objects are equivalent if and only if there is a finite string of
arrows connecting them in the ambient category, and we show that there
are functions
\begin{equation*}
\Phi: \pi_{0}(\mathbf{Gerbe}(\mathcal{C})) \leftrightarrows 
\pi_{0}H(\ast,\mathbf{Grp}(\mathcal{C})): \Psi
\end{equation*}
which are inverse to each other. Here, $\Phi$ is induced by a canonical
cocycle construction which is introduced in Example \ref{ex 12}, and
$\Psi$ is defined by a generalized Grothendieck construction, which is
the subject of much of Section 2.

The $2$-groupoid diagram $\mathbf{Grp}(\mathcal{C})$ has subobjects
which are honest presheaves of $2$-groupoids. Examples include the
sheaf of $2$-groupoids $G_{\ast}$ associated to a sheaf of groups of
$G$ on $\mathcal{C}$: it has one object, the sheaf of $1$-cells is the
sheaf of automorphisms of $G$, and its sheaf of $2$-cells is the sheaf
of homotopies (or conjugations) of automorphisms. More generally, any
presheaf of sheaves of groups in $\mathbf{Grp}(\mathcal{C})$
determines a full subobject $\mathcal{F} \subset
\mathbf{Grp}(\mathcal{C})$ which is a presheaf of $2$-groupoids, and
one can discuss the homotopy type of $\mathcal{F}$ and its classifying
object $B\mathcal{F}$ in simplicial presheaves. It is shown in Theorem
\ref{thm 23}
that there is a one to one correspondences
\begin{equation*}
\pi_{0}H(\ast,\mathcal{F}) \cong
\pi_{0}(\mathcal{F}-\mathbf{Gerbe})
\end{equation*}
between the set of path components of $2$-cocycles taking values in
the presheaf of $2$-groupoids $\mathcal{F}$ and path components of the
category $\mathcal{F}-\mathbf{Gerbe}$ of gerbes locally equivalent to
sheaves of groups appearing in $\mathcal{F}$. By the result
relating path components of cocycle categories to morphisms in the
homotopy category displayed above, both of these objects are then in
bijective correspondence with the set $[\ast,B\mathcal{F}]$ of
morphisms in the homotopy category of simplicial presheaves --- this
statement appears formally as Corollary \ref{cor 24}. The bijection of
path components in the statement of Theorem \ref{thm 23} is a
restriction of the bijection of Theorem \ref{thm 20}.

In the special case where $\mathcal{F} = G_{\ast}$ for some sheaf of
groups $G$, Theorem \ref{thm 23} says that gerbes locally equivalent
to $G$ are classified up to weak equivalence by morphisms
$[\ast,BG_{\ast}]$ in the homotopy category of simplicial
presheaves. This result was originally proved, in a very different
form, by Breen \cite{Br}.

Finally, the presheaf of $2$-groupoids has a fundamental groupoid $\pi
\mathcal{F}$ and a canonical morphism $\mathcal{F} \to \pi
\mathcal{F}$. In the case where $\mathcal{F}=G_{\ast}$, the
fundamental groupoid $\pi G_{\ast}$ is the sheaf of outer
automorphisms of $G$. The fundamental groupoid $\pi \mathcal{F}$ has a
functorial stack completion $\pi \mathcal{F} \to \St(\pi
\mathcal{F})$, and $\St(\pi \mathcal{F})$ is the stack of bands
(liens) for $\mathcal{F}$. Suppose that the band $L$ is a fixed choice
of global section of $\St(\pi \mathcal{F})$, and consider the homotopy
fibre $F_{L}$ of the composite
\begin{equation*}
B\mathcal{F} \to B\pi \mathcal{F} \to B\St(\pi\mathcal{F}).
\end{equation*}
Theorem \ref{thm 27} (see also Corollary \ref{cor 28}) identifies the
set of morphisms $[\ast,BF_{L}]$ in the homotopy category with path
components in a suitably defined category of $L$-gerbes. In other
words, Giraud's non-abelian invariant $H^{2}(\mathcal{C},L)$ is
isomorphic to $[\ast,BF_{L}]$. Once again, the real thrust of the
proof is to identify the set of path components of the cocycle
category $H(\ast,F_{L})$ with path components in $L$-gerbes, and then
use the general result about cocycles to conclude that $[\ast,F_{L}]
\cong \pi_{0}H(\ast,F_{L})$.

Theorem \ref{thm 20}, Theorem \ref{thm 23} and Theorm \ref{thm 27} are
the main results of this paper. The demonstration of these results
appear in Sections 3 and 4, but depend on some generalities about
groupoids enriched in simplicial sets and presheaves of $2$-groupoids
which appear in Section 1, as well as the discussion of the generalized
Grothendieck construction of Section 2.  

\tableofcontents

\section{Simplicial groupoids}

There are various equivalent ways to define a groupoid enriched in simplicial
sets. We shall initially take the point of view that such an object
$H$ is a simplicial groupoid such that the simplicial set $\Ob(H)$ of
objects is simplicially discrete, or just a set. The morphisms
$\Mor(H)$ is a simplicial sets and the source, target $s,t: \Mor(H)
\to \Ob(H)$ and identity maps $e: \Ob(H) \to \Mor(H)$ are all
simplicial set maps. The notation $H_{n}$ will refer to the associated
groupoid in simplicial degree $n$.

For objects $x,y$ of $H$ the simplicial set of morphisms
$H(x,y)$ is defined by the pullback diagram
\begin{equation*}
\xymatrix{
H(x,y) \ar[r] \ar[d] & \Mor(H) \ar[d]^{(s,t)} \\
\ast \ar[r]_-{(x,y)} & \Ob(H) \times \Ob(H)
}
\end{equation*}
where $\ast = \Delta^{0}$ defines the one-point (terminal) simplicial
set.

A $2$-groupoid is a groupoid enriched in groupoids. Equivalently, a
simplicial groupoid $G$ is a groupoid enriched in simplicial sets such
that the simplicial set $\Mor(G)$ is the nerve of a groupoid.

We shall routinely write $BH$ for both the bisimplicial classifying
space $n \mapsto BH_{n}$ associated to $H$ and its associated diagonal
simplicial set $dBH$. The vertical simplicial presheaf $BH_{n}$ in
horizontal degree $n$ is the iterated fibre product
\begin{equation*}
\Mor(H) \times_{t,s} \Mor(H) \times_{t,s} \dots \times_{t,s} \Mor(H),
\end{equation*} 
This simplicial set models composable strings of morphisms of
length $n$, and is the inverse limit for a diagram
\begin{equation*}
\xymatrix@C=8pt{
\Mor(H) \ar[dr]^{t} && \Mor(H) \ar[dl]_{s} \ar[dr]^{t} 
&& \Mor(H) \ar[dl]_{s} & \dots \\
& \Ob(H) && \Ob(H)
}
\end{equation*}
involving $n$ copies of the morphism object $\Mor(H)$.

The vertices of $BH$ are the objects of $H$, and two vertices of $BH$
are in the same path component if and only if they are in the same
path component of the space $BH_{1}$, or in the same path component of
the groupoid $H_{1}$ in simplicial degree $1$. It is well known and
easily seen that the degeneracy morphism $H_{0} \to H_{n}$ induces a
bijection $\pi_{0}BH_{0} \cong \pi_{0}BH_{n}$ for all $n \geq 0$. It
follows that there are natural isomorphisms
\begin{equation}\label{eq 1}
\pi_{0}BH \cong \pi_{0}BH_{1} \cong \pi_{0}BH_{n}
\end{equation}
for all $n \geq 0$. We shall say that $H$ is {\it connected} if $BH$
is a path-connected simplicial set. More generally, $\pi_{0}H$ will
often be written to denote $\pi_{0}BH$, and will be called the set of
{\it path components} of $H$.

There
is a simplicial groupoid $H^{\mathbf{1}}$ whose objects in simplicial
degree $n$ are morphisms $h: x \to y$ of $H_{n}$ and whose morphisms
$h \to h'$ are the commutative squares
\begin{equation*}
\xymatrix{
x \ar[r]^{h} \ar[d]_{\alpha} & y \ar[d]^{\beta} \\
x' \ar[r]_{h'} & y'
}
\end{equation*}
in $H_{n}$. There is a simplicial groupoid functor $(s,t):
H^{\mathbf{1}} \to H \times H$ which is defined in degree $n$ by
sending the square diagram above to the pair of morphisms $(x
\xrightarrow{\alpha} x',y \xrightarrow{\beta} y')$. 
The
two projections $s,t: H^{\mathbf{1}} \to H$ are weak equivalences,
because they are weak equivalences in each simplicial degree.

\begin{lemma}\label{lem 1}
Suppose that $H$ is a groupoid enriched in simplicial
sets. Then the pullback diagram
\begin{equation*}
\xymatrix{
\Mor(H) \ar[r] \ar[d] & BH^{\mathbf{1}} \ar[d] \\
\Ob(H) \times \Ob(H) \ar[r] & BH \times BH
}
\end{equation*}
is homotopy cartesian.
\end{lemma}

\begin{proof}
Suppose first that $H$ is an ordinary groupoid. There is a functor $\hom: H
\times H \to \mathbf{Set}$ defined by $(x,y) \mapsto H(x,y)$ and which
sends the pair of morphisms $(\alpha: x \to x',\beta: y \to y')$ to
the function
\begin{equation*}
H(x,y) \to H(x',y')
\end{equation*}
defined by sending $f:x \to y$ to the composite
\begin{equation*}
x' \xrightarrow{\alpha^{-1}} x \xrightarrow{f} y \xrightarrow{\beta} y'.
\end{equation*}
Observe that there is an isomorphism
\begin{equation*}
\hocolim_{H \times H} \hom \cong B(H^{\mathbf{1}}).
\end{equation*}

Suppose now that $H$ is a groupoid enriched in simplicial sets.
Applying the construction of the previous paragraph in all simplicial
degrees gives a simplicial functor $\hom: H \times H \to
s\mathbf{Set}$ and an isomorphism
\begin{equation*}
\hocolim_{H \times H}\hom \cong B(H^{\mathbf{1}}).
\end{equation*}
There is a homotopy cartesian diagram
\begin{equation*}
\xymatrix{
X \ar[r] \ar[d] & \hocolim_{G}X \ar[d] \\
\Ob(G) \ar[r] & BG
}
\end{equation*}
for all simplicial set-valued diagrams $X$ defined on groupoids $G$
enriched in simplicial sets \cite[Lemma 2]{J2}, \cite{Mo}, and this
specializes to give the homotopy cartesian diagram required by the
statement of the Lemma.
\end{proof}

\begin{lemma}\label{lem 2}
If $f: G \to H$ is a morphism of groupoids enriched in simplicial sets
and if $X: H \to s\mathbf{Set}$ is a simplicial diagram defined on
$H$, then the induced diagram
\begin{equation*}
\xymatrix{
\hocolim_{G} X\cdot f \ar[r] \ar[d] & \hocolim_{H} X \ar[d] \\
BG \ar[r] & BH
}
\end{equation*}
is homotopy cartesian.
\end{lemma}

\begin{proof}
Suppose given a factorization
\begin{equation*}
\xymatrix{
\hocolim_{H}X \ar[r]^-{j} \ar[dr] & Z \ar[d]^{p} \\
& BH
}
\end{equation*}
of the canonical map $\hocolim_{H} X \to BH$ such that $j$ is a weak
equivalence and $p$ is a fibration. The squares in the picture
\begin{equation*}
\xymatrix{
X \cdot f \ar[r] \ar[d] & X \ar[r] \ar[d] & \hocolim_{H} X \ar[d] \\
\Ob(G) \ar[r] & \Ob(H) \ar[r] & BH
}
\end{equation*}
are homotopy cartesian, so that the induced map
\begin{equation*}
\xymatrix@C=8pt{
X \cdot f \ar[rr] \ar[dr] && \Ob(G) \times_{BH} Z \ar[dl] \\
& \Ob(G)
}
\end{equation*}
is a weak equivalence of objects over $\Ob(G)$. It follows that the induced map
\begin{equation*}
\xymatrix@C=8pt{
\hocolim_{G} X \cdot f \ar[rr] \ar[dr] && BG \times_{BH} Z \ar[dl] \\
& BG
}
\end{equation*}
induces a weak equivalence on all homotopy fibres. This map is also a
homotopy colimit of a comparison of diagrams made up of the homotopy
fibres of the respective maps, and is therefore a weak equivalence.
\end{proof}

\begin{corollary}\label{cor 3}
Suppose that $H$ is a groupoid enriched in simplicial
sets. Then the square
\begin{equation*} 
\xymatrix{ 
H(x,y) \ar[r] \ar[d] &
B(H/y) \ar[d] \\ 
\ast \ar[r]_{x} & BH 
}
\end{equation*}
is homotopy cartesian.
\end{corollary}

\begin{proof}
In the diagram of pullback squares
\begin{equation*}
\xymatrix{
H(x,y) \ar[r] \ar[d] & B(H/y) \ar[r] \ar[d] & B(H^{\mathbf{1}}) \ar[d] \\
\ast \ar[r]_{x} & BH \ar[r]_-{(1,y)} & BH \times BH
}
\end{equation*}
the square on the right is homotopy cartesian, on account of Lemma
\ref{lem 2} applied to the composite functor
\begin{equation*}
H \xrightarrow{(1,y)} H \times H \xrightarrow{hom} s\mathbf{Set}.
\end{equation*} 
Lemma \ref{lem 1} implies that the composite square is homotopy
cartesian, and the desired result follows.
\end{proof}

Let $\pi H$ denote the groupoid of path components of a groupoid $H$
enriched in simplicial sets. The object $\pi H$ will typically be
called the {\it path component} groupoid of $H$. It has the same objects as
$H$, and the set of morphisms from $x$ to $y$ is the set $\pi H(x,y)$
of path components of the simplicial set $H(x,y)$. There is a
canonical map $\eta: H \to \pi H$ which is the identity on objects and
is the canonical map $\Mor(H) \to \pi_{0}\Mor(H)$ on morphisms. The
morphism $\eta$ is one of the canonical maps for an adjunction: the
functor $H \mapsto \pi H$ is left adjoint to the inclusion of
groupoids in simplicial groupoids.

\begin{corollary}
The induced map $\eta: BH \to B\pi H$ induces an isomorphism on path
components and all fundamental groups, so that $\pi H$ is naturally
weakly equivalent to the fundamental groupoid of $BH$.
\end{corollary}

\begin{proof}
The morphisms $H(x,x) \to \pi H(x,x)$ induce isomorphisms in path
components, and so the map $BH \to B\pi H$ induces isomorphisms in
path components of all loop spaces, by Corollary \ref{cor 3}. It
follows that all homomorphisms $\pi_{1}(BH,x) \to \pi_{1}(B\pi H,x)$
are isomorphisms. The claim that $\pi_{0}BH \to \pi_{0}B\pi H$ is a
bijection follows from (\ref{eq 1}) and the observation that the
function $\pi_{0}BH_{0} \to \pi_{0}B\pi H$ is a bijection.
\end{proof}

Suppose now that $\mathcal{C}$ is a small Grothendieck site.

If $G$ is a presheaf of groupoids on $\mathcal{C}$ and $x,y$ are
objects of $G(U)$, there is a presheaf $G(x,y)$ of homomorphisms from
$x$ to $y$ on $\mathcal{C}/U$. Write $G_{x} = G(x,x)$ for the presheaf
of automorphisms of $x$ in $G$, and let $\tilde{G}_{x}$ denote the
associated sheaf of automorphisms on $\mathcal{C}/U$.

Say that a presheaf of groupoids $G$ is a {\it \v{C}ech object} if the
canonical map $G \to \pi_{0}G$ is a local weak equivalence, where
$\pi_{0}G = \pi_{0}BG$ is the presheaf of path components of $G$.  

In particular, an ordinary groupoid $H$ is a \v{C}ech groupoid if the
groupoid moprhism $H \to \pi_{0}H$ is a weak
equivalence. Equivalently, $H$ is a \v{C}ech groupoid if and only if
there is at most one morphism between any two objects of $H$.

\begin{example}
The \v{C}ech groupoid $C(p)$ for a function $p: X \to Y$
has the objects $\Ob(C(p)) = X$, and there is a morphism $x \to y$ in
$C(p)$ if and only if $p(x) = p(y)$ in $Y$. There is a canonical
bijection $\pi_{0}C(p) \cong p(X)$.

This construction is natural, and therefore applies to morphisms $p: X
\to Y$ of presheaves on a site. If $p$ is an epimorphism of sheaves,
the simplicial presheaf map $BC(p) \to Y$ is the \v{C}ech resolution
of $Y$ corresponding to the epimorphism $p$, and is a local weak
equivalence. This construction specializes to the standard \v{C}ech
resolution when applied to an epimorphism
$p: \bigsqcup_{\alpha} U_{\alpha} \to Y$ arising from a covering.
\end{example}

\begin{lemma}\label{lem 6}
Suppose that $H$ is a presheaf of $2$-groupoids, and let $\eta: H \to
\pi H$ be the canonical map to the presheaf of path component
groupoids. Then $\eta$ is a local weak equivalence if and
only if all presheaves of groupoids $H(x,y)$ are \v{C}ech objects.
\end{lemma}

\begin{proof}
If the map $\eta: H \to \pi H$ is a local weak equivalence, then the map
\begin{equation*}
\Mor(H) \to \pi_{0}\Mor(H) 
\end{equation*}
is a local weak equivalence over $\Ob(H)
\times \Ob(H)$, by Lemma \ref{lem 1}. The object $\Ob(H) \times
\Ob(H)$ is a constant simplicial presheaf, so pullback along any map
$Z \to \Ob(H) \times \Ob(H)$ preserves local weak equivalences over
$\Ob(H) \times \Ob(H)$ \cite{J2}. In particular, for all choices $x,y
\in \Ob(U)(U)$, the induced map $H(x,y) \to \pi_{0}H(x,y)$ is a local
weak equivalence of presheaves of groupoids over $\mathcal{C}/U$. 

In general, one can show that a map
\begin{equation*}
\xymatrix@C=8pt{
Z \ar[rr]^{f} \ar[dr] && W \ar[dl] \\
& A
}
\end{equation*}
of simplicial presheaves fibred over a presheaf $A$ is a local weak
equivalence if and only if it induces weak equivalences $Z_{x} \to
W_{x}$ of simplicial presheaves on $\mathcal{C}/U$ for all $x \in
A(U), U \in \mathcal{C}$. One implication involved in this statement
we already know about, from \cite{J2}. For the other, if we know that
all induced maps on fibres are local weak equivalences, we can replace
$f$ by a sectionwise Kan fibration, and check local lifting with
respect to all inclusions $\partial\Delta^{n} \subset \Delta^{n}$,
which must take place in individual fibres.

Thus, if all morphism groupoids $H(x,y)$ are \v{C}ech objects, then
the map $\Mor(H) \to \pi_{0}\Mor(H)$ is a local weak equivalence of
simplicial presheaves over the presheaf
$\Ob(H) \times \Ob(H).$ 
It follows that all maps
\begin{equation*}
\begin{aligned}
\Mor(H) &\times_{t,s} \Mor(H) \times_{t,s} \dots \times_{t,s} \Mor(H) \\
&\to
\pi_{0}\Mor(H) \times_{t,s} \pi_{0}\Mor(H) \times_{t,s} \dots \times_{t,s} \pi_{0}\Mor(H) 
\end{aligned}
\end{equation*}
of iterated fibre products over $\Ob(H)$ are local weak
equivalences. These are the comparison maps of vertical simplicial
presheaves making up the comparison $BH \to B\pi H$ of bisimplicial
presheaves, and one concludes that this map is a
local weak equivalence.
\end{proof}

\begin{lemma}\label{lem 7}
A presheaf of groupoids $H$ is a \v{C}ech object if and only if for
every two morphisms $f,g: x \to y$ in $H(U)$ there is a covering sieve
$R \subset \hom(\ ,U)$ such that $\phi^{\ast}f = \phi^{\ast}g$ for all
$\phi: V \to U$ in $R$.
\end{lemma}

\begin{proof}
Suppose that $H \to \tilde{H}$ is the canonical map taking values in
the associated sheaf of groupoids $\tilde{H}$. Then $H$ is a \v{C}ech object
if and only if all sheaves $\tilde{H}(x,x)$ of automorphisms of
$\tilde{H}$ are trivial in the sense that the canonical sheaf map
$\tilde{H}(x,x) \to \ast$ are isomorphisms. This is equivalent to the
assertion that all presheaf maps $H(x,x) \to \ast$ are local
monomorphisms. 

Thus, suppose that $H$ is a \v{C}ech object, and suppose given $f,g: x
\to y$ in $H(U)$ the composite $g^{-1}f \in H(x,x)(U)$, and there is a
covering sieve $R \subset \hom(\ ,U)$ such that $\phi^{\ast}(g^{-1}f)
= 1_{\phi^{\ast}x}$ for all $\phi: V \to U$ in $R$. But then
$\phi^{\ast}(g) = \phi^{\ast}(f)$ for all $\phi \in R$.

The converse is clear: the local coincidence of all $f,g: x \to y$
means that all presheaf maps $H(x,y) \to \ast$ are local
monomorphisms, and so all sheaf maps $\tilde{H}(x,x) \to \ast$ are
isomorphisms.
\end{proof}

\begin{lemma}\label{lem 8}
Suppose that $A$ is a presheaf of $2$-groupoids, and that
$\pi_{0}A$ is its presheaf of path components. Then the canonical
map $A \to \pi_{0}A$ is a local weak equivalence if and only
if all presheaves of groupoids $A(x,y)$ and the path component
groupoid $\pi A$ are \v{C}ech objects.
\end{lemma}

\begin{proof}
Suppose that $A \to \pi_{0}A$ is a local weak
equivalence. Then all sheaves of homotopy groups for $BA$ are trivial,
and so Lemma \ref{lem 1} implies that all maps $A(x,y) \to \ast$ are
local weak equivalences. In particular, all $A(x,y)$ are \v{C}ech
objects. But then $\eta: A \to \pi A$ is a local weak equivalence by
Lemma \ref{lem 6}, and so the induced map $\pi A \to \pi_{0}(\pi A)$
is a weak equivalence, so that the
presheaf of groupoids $\pi A$ is a \v{C}ech object.

Suppose conversely that all $A(x,y)$ and $\pi A$ are \v{C}ech
objects. Then Lemma \ref{lem 6} implies that $A \to \pi A$ is a
local weak equivalence, and then the map
\begin{equation*}
\pi A \to \pi_{0}(\pi A) \cong \pi_{0}A
\end{equation*}
is a local weak equivalence. It follows that the map $A \to
\pi_{0}A$ is a composite of two local weak equivalences.
\end{proof}

\section{The Grothendieck construction}

Let $\mathbf{cat}$ denote the $2$-category whose $0$-cells are the
small categories, whose $1$-cells are the functors between small
categories, and whose $2$-cells are the homotopies of functors.

Suppose given a $2$-category morphism $F: A \to \mathbf{cat}$, such
that $A$ is a small category enriched in groupoids. This morphism has
an associated ``Grothendieck construction'', which is a category
$E_{A}F$ that is constructed as follows.

Consider the collection of pairs $(x,i)$ where $i$ is an object or
$0$-cell of $A$ and $x \in F(i)$. Look at all pairs
\begin{equation*}
(f,\alpha): (x,i) \to (y,j)
\end{equation*}
where $\alpha: i \to j$ is a $1$-cell of $A$ and $f: \alpha_{\ast}(x)
\to y$ is a morphism of $F(j)$. Say that two such pairs
\begin{equation*}
(f,\alpha), (f',\alpha'): (x,i) \to (y,j)
\end{equation*}
are equivalent if there is a $2$-cell $h:\alpha \to \alpha'$ in $A$
such that the diagram
\begin{equation*}
\xymatrix@R=6pt{
\alpha_{\ast}(x) \ar[dr]^{f} \ar[dd]_{F(h)} & \\
& y \\
\alpha'_{\ast}(x) \ar[ur]_{f'}
}
\end{equation*}
commutes, where $F(h)$ is the homotopy associated to the $2$-cell $h$
by $F$. This is an equivalence relation, since the homotopies $h$
are isomorphisms in the groupoids $A(x,y)$.
Write $[(f,\alpha)]$ for the equivalence class containing the pair
$(f,\alpha)$.

Suppose given strings
\begin{equation*}
(x,i) \xrightarrow{(f,\alpha)} (y,j) \xrightarrow{(g,\beta)} (z,k)
\end{equation*}
and
\begin{equation*}
(x,i) \xrightarrow{(f',\alpha')} (y,j) \xrightarrow{(g',\beta')} (z,k)
\end{equation*}
and suppose that $(f,\alpha) \overset{h_{1}}{\simeq} (f',\alpha')$ and
$(g,\beta) \overset{h_{2}}{\simeq} (g',\beta')$ via the displayed
homotopies. Then there is a commutative diagram
\begin{equation*}
\xymatrix{
\beta_{\ast}\alpha_{\ast}(x) \ar[r]^{\beta_{\ast}(f)} \ar[d]_{F(h_{2})} 
& \beta_{\ast}(y) \ar[r]^{g} \ar[d]_{F(h_{2})} & z \\
\beta'_{\ast}\alpha_{\ast}(x) \ar[r]^{\beta'_{\ast}(f)} 
\ar[d]_{\beta'_{\ast}(F(h_{1}))} 
& \beta'_{\ast}(y) \ar[ur]_{g'} & \\
\beta'_{\ast}\alpha'_{\ast}(x) \ar[ur]_{\beta'_{\ast}(f')}
}
\end{equation*}
The composite homotopy
$\beta'_{\ast}(F(h_{1}))(F(h_{2})\alpha_{\ast})$ is the image of the
composite $2$-cell $h_{2} \ast h_{1}$ under the morphism $F$. It
follows that the assignment
\begin{equation*}
[(g,\beta)]\cdot [(f,\alpha)] = [(g\beta_{\ast}(f),\beta\alpha)]
\end{equation*}
gives a well defined law of composition. This composition law is
associative, and has $2$-sided identities. Write $E_{A}F$ for the
corresponding category.

\begin{remark}
This category $E_{A}F$ should seem familiar. Suppose that $I$ is
a small category, and let $x,y$ be objects of $I$. There is a small
category $I_{s}(x,y)$ whose objects are the functors $\theta:
\mathbf{n} \to I$ (strings of length $n$) such that $\theta(0) = x$
and $\theta(n) = y$. A morphism $\theta \to \gamma$ of $I_{s}(x,y)$ is
a commutative diagram of functors
\begin{equation*}
\xymatrix@R=8pt{ \mathbf{n} \ar[dr]^{\theta} \ar[dd]_{\alpha} & \\ 
& I \\ 
\mathbf{m} \ar[ur]_{\gamma}
}
\end{equation*}
such that the ordinal number map $\alpha$ is end-point preserving in
the sense that $\alpha(0) = 0$ and $\alpha(n) = m$. Concatenation of
strings defines a composition law $I_{s}(x,y) \times I_{s}(y,z) \to
I_{s}(x,z)$, and so there is a $2$-category $I_{s}$ with the same
objects as $I$ and a canonical weak equivalence $I_{s} \to I$ (see
also \cite[IX.3.2]{GJ}). Write $GI_{s}$ for the category enriched in
groupoids, having the same $0$-cells as $I$, and such that the
groupoid $GI_{s}(x,y)$ is the free groupoid on the category
$I_{s}(x,y)$. The groupoid $GI_{s}(x,y)$ is a \v{C}ech groupoid with
path components given by the set $I(x,y)$ of morphisms from $x$ to $y$
in $I$.

A pseudo-functor $F$ defined on $I$ and taking values in small
categories can be identified with a $2$-category morphism $F: GI_{s}
\to \mathbf{cat}$ \cite[IX.3.3]{GJ}, and one can show that the
Grothendieck construction $E_{F}GI_{s}$ as defined above is isomorphic
to the standard Grothendieck construction for the
pseudo-functor $F$. 

The Grothendieck construction $E_{A}F$ given here for $2$-catgory
morphisms defined on categories $A$ enriched in groupoids generalizes
the usual construction for pseudo-functors, but
there appears to be no corresponding construction for lax functors.
\end{remark}

We shall henceforth specialize to $2$-category morphisms $F: A \to
\mathbf{cat}$ which are defined on small $2$-groupoids $A$.
In this case, there is a canonical functor $p: E_{A}F \to \pi A$, which
is defined by the assignment $(x,i) \mapsto i$.

\begin{lemma}\label{lem 10}
Suppose that $F: A \to \mathbf{cat}$ is a $2$-category morphism, where
$A$ is a $2$-groupoid.  Suppose that $[(f,\alpha)]: (x,i) \to (y,j)$
is a morphism of $E_{A}F$ such that $f: \alpha_{\ast}(x) \to y$ is an
invertible morphism of $F(j)$. Then $[(f,\alpha)]$ is invertible in
$E_{F}A$.
\end{lemma}

\begin{proof}
The inverse of $[(f,\alpha)]$ is represented by
$[(\alpha^{-1}_{\ast}(f^{-1}),\alpha^{-1})]$.
\end{proof}

\begin{corollary}
If the $2$-category morphism $F: A \to \mathbf{cat}$ takes values in
groupoids, then $E_{A}F$ is a groupoid.
\end{corollary}

A diagram
\begin{equation*}
B  \underset{\simeq}{\xleftarrow{\alpha}} A \xrightarrow{F} \mathbf{cat}
\end{equation*}
such that $\alpha: A \to B$ is a weak equivalence of $2$-groupoids is a
{\it $2$-cocycle} taking values in small categories. A
morphism of $2$-cocycles is a commutative diagram of functors
\begin{equation*}
\xymatrix@R=8pt{
& A \ar[dl]_{\alpha}^{\simeq} \ar[dd] \ar[dr]^{F} & \\
B && \mathbf{cat} \\
& A' \ar[ul]^{\alpha'}_{\simeq} \ar[ur]_{F'}
}
\end{equation*}
and the corresponding $2$-cocycle category is denoted by
$H(B,\mathbf{cat})$. There are analogous definitions for $2$-cocycles
and $2$-cocycle categories taking values in 
of groups and small groupoids. These $2$-cocycle categories 
are typically not small.

In particular, write $\mathbf{grp}$ for the $2$-groupoid whose objects
are the groups, whose $1$-cells are the isomorphisms of groups $G \to
H$, and whose $2$-cells are the homotopies of isomorphisms, and
suppose now that there is a $2$-cocycle
\begin{equation*}
\pi A \underset{\simeq}{\xleftarrow{\eta}} A \xrightarrow{K} \mathbf{grp}
\end{equation*}
taking values in groups. Then the associated Grothendieck
construction $E_{A}K$ can be identified with a category having as
objects all $i \in \Ob(A)$ and with morphisms consisting of
equivalence classes of pairs $(f,\alpha): i \to j$, where $\alpha: i
\to j$ is a $1$-cell of $A$ and $f \in K(j)$. In this case, there is a
relation $(f,\alpha) \sim (g,\beta)$ if
the diagram
\begin{equation*}
\xymatrix@R=6pt{
\ast \ar[dr]^{f} \ar[dd]_{h_{\ast}} & \\
& \ast \\
\ast \ar[ur]_{g} &
}
\end{equation*}
commutes in the group $K(j)$, where conjugation by $h_{\ast}$
defines the image of the unique $2$-cell $\alpha \to \beta$.

The category $E_{A}K$ is a groupoid by Lemma \ref{lem 10}. 

\begin{example}\label{ex 12}
Suppose that $G$ is a groupoid. The {\it resolution $2$-groupoid}
$R(G)$ has the same objects and $1$-cells as $G$, and has a unique
$2$-cell $f \to g$ between any two morphisms $f,g: x \to y$ of $G$.
The path component groupoid $\pi R(G)$ of $R(G)$ is a \v{C}ech groupoid,
and the natural maps
\begin{equation*}
BR(G) \to B\pi R(G) \to \pi_{0}B(\pi R(G))
\end{equation*} 
weak equivalences. There are natural
bijections
\begin{equation*}
\pi_{0}G = \pi_{0}BG \cong \pi_{0}B(\pi R(G)).
\end{equation*}

There is a canonical morphism $F(G): R(G) \to \mathbf{grp}$ which
takes the object $x \in G(U)$ to the group $G_{x} = G(x,x)$ on
$\mathcal{C}/U$, takes a $1$-cell $f: x \to y$ to the isomorphism
$G_{x} \to G_{y}$ which is defined by conjugation by $f$, and takes
the $2$-cell $f \to g$ to the homotopy defined by
conjugation by the element $gf^{-1} \in G_{y}$.
It follows that $G$ determines a canonical $2$-cocycle
\begin{equation*}
\pi_{0}G \xleftarrow{\simeq} R(G) \xrightarrow{F(G)} \mathbf{grp}.
\end{equation*}
\end{example}

\begin{lemma}\label{lem 13}
There is a natural isomorphism of groupoids $\psi: E_{R(G)}F(G)
\xrightarrow{\cong} G$ which is defined fibrewise over
$\pi R(G)$ in the sense that there is a commutative diagram
\begin{equation*}
\xymatrix@C=6pt{
E_{R(G)}F(G) \ar[rr]^-{\psi}_-{\cong} \ar[dr] && G \ar[dl] \\
& \pi R(G)
}
\end{equation*}
\end{lemma}

\begin{proof}
The functor $\psi$ is the identity on objects. It is defined on
morphisms by sending the pair $(f,\alpha)$ to the composite $f\cdot
\alpha$ in $G$. If $(f,\alpha) \sim (g,\beta)$ and $\alpha \to \beta$ is the unique $2$-cell in $R(G)$,
then the diagram
\begin{equation*}
\xymatrix@R=6pt{
& j \ar[dd]^{\beta\alpha^{-1}} \ar[dr]^{f} & \\
i \ar[ur]^{\alpha} \ar[dr]_{\beta} && j \\
& j \ar[ur]_{g} &
}
\end{equation*}
commutes in $G$, so that $f\cdot \alpha = g\cdot \beta$ and the
assignment $[(f,\alpha)] \mapsto f\cdot \alpha$ is well defined. The
assignment is functorial, because the diagram
\begin{equation*}
\xymatrix@R=6pt{
&& k \ar[dr]^{\beta_{\ast}(f)} && \\
& j \ar[ur]^{\beta} \ar[dr]^{f} && k \ar[dr]^{g} & \\
i \ar[ur]^{\alpha} && j \ar[ur]^{\beta} && k
}
\end{equation*}
commutes in $G$. 

The functor $\psi$ plainly induces surjective functions
\begin{equation*}
\psi: \hom_{E_{R(G)}F(G)}(i,j) \to \hom_{G}(i,j)
\end{equation*}
Finally, if the diagram
\begin{equation*}
\xymatrix@R=6pt{
& j \ar[dr]^{f} & \\
i \ar[ur]^{\alpha} \ar[dr]_{\beta} && j \\
& j \ar[ur]_{g} &
}
\end{equation*}
commutes in the groupoid $G$ then $g \cdot(\beta\alpha^{-1}) = f$, so
that $(f,\alpha) \sim (g,\beta)$, and $\psi$ is injective on
morphisms.
\end{proof}

\section{Cocycle classification of gerbes}

A {\it gerbe} $G$ is a locally connected presheaf of groupoids. A
{\it morphism of gerbes} is a local weak equivalence $G \to H$ of
presheaves of groupoids. We shall write $\mathbf{Gerbe}(\mathcal{C})$
for the category of gerbes and morphisms of gerbes.

\begin{remark}
If $G$ is a gerbe and $x$ is a global section of $\Ob(G)$, then the
inclusion map $G_{x} \to G$ is a local weak equivalence. It follows
that every gerbe $H$ is locally equivalent to a presheaf of groups, in
the sense that there is a covering $U \to \ast$ by objects $U \in
\mathcal{C}$ and section $x_{U} \in \Ob(H)(U)$ such that the morphisms
$H_{x_{U}} \to H\vert_{U}$ are local weak equivalences over
$\mathcal{C}/U$ for all $U$ in the covering.
\end{remark}

\begin{remark}
Suppose that $E$ is a presheaf, and identify $E$ with a presheaf of
discrete groupoids. An {\it $E$-gerbe} is a morphism $G \to E$ of
presheaves of groupoids such that the associated presheaf map
$\pi_{0}G \to E$ induces an isomorphism $\tilde{\pi}_{0}G \cong
\tilde{E}$ of associated sheaves.  A morphism of $E$-gerbes is a
commutative diagram
\begin{equation*}
\xymatrix@C=6pt{
G \ar[rr]^{f} \ar[dr] && H \ar[dl] \\
& E
}
\end{equation*}
such that the morphism $f: G \to H$ is a local weak equivalence of
presheaves of groupoids. Write $\mathbf{Gerbe}_{E}(\mathcal{C})$ for
the corresponding category.  Categories of $E$-gerbes do appear in
applications --- see \cite[p.22]{LM}.

There is an equivalence of categories 
\begin{equation*}
\mathbf{Gerbe}_{E}(\mathcal{C}) \simeq \mathbf{Gerbe}(\mathcal{C}/E)
\end{equation*}
between the category
of $E$-gerbes on $\mathcal{C}$ and the
category of gerbes for the fibred site
$\mathcal{C}/E$. Equivalences of this sort are discussed at length in
\cite{J1.5}.
Classification results for
$E$-gerbes can therefore be deduced from classification results for gerbes on
the site $\mathcal{C}/E$.
\end{remark}

We shall write $\mathbf{Grp}(\mathcal{C})$ for the following monster:
it is a contravariant diagram defined on $\mathcal{C}$ and taking
values in $2$-groupoids, such that the $0$-cells of
$\mathbf{Grp}(\mathcal{C})(U)$ are the sheaves of groups on
$\mathcal{C}/U$, the $1$-cells are the isomorphisms of sheaves of
groups on $\mathcal{C}/U$, and the $2$-cells are the (global)
homotopies of sheaf isomorphisms. $\mathbf{Grp}(\mathcal{C})$ is not a
presheaf of groupoids, because it does not take values in small
groupoids.

If $G$ is a gerbe, then the corresponding resolution $2$-groupoid $R(G)$
(Example \ref{ex 12}) is weakly equivalent to a point in the sense
that the map $R(G) \to \ast$ is a local weak equivalence of presheaves
of $2$-groupoids.  There is a canonical morphism $F(G): R(G) \to
\mathbf{Grp}(\mathcal{C})$ for which the $0$-cell $x \in R(G)(U)$ is
mapped to the sheaf of groups $\tilde{G}_{x}$, the $1$-cell $\alpha: x
\to y$ is mapped to the sheaf isomorphism $c_{\alpha}: \tilde{G}_{x}
\to \tilde{G}_{y}$ on $\mathcal{C}/U$ which is defined by conjugation
by the global section $\alpha$, and each $2$-cell $h: \alpha \to
\beta$ of $1$-cells $x \to y$ maps to conjugation by the image of $h
\in \tilde{G}_{y}(U)$ in global sections of $\tilde{G}_{y}$. In this
way, each gerbe $G$ has a canonically associated $2$-cocycle
\begin{equation*}
\ast \xleftarrow{\simeq} R(G) \xrightarrow{F(G)} \mathbf{Grp}(\mathcal{C}).
\end{equation*}

Write $H(\ast,\mathbf{Grp}(\mathcal{C}))$ for the category of
$2$-cocycles taking values in the $2$-groupoid object
$\mathbf{Grp}(\mathcal{C})$.

The assignment of the cocycle $F(G): R(G) \to
\mathbf{Grp}(\mathcal{C})$ to the gerbe $G$ is not
functorial. It is true, however, that a map $G \to H$ of
gerbes induces a $2$-groupoid morphisms $f_{\ast}: R(G) \to R(H)$.
The sheaf isomorphisms $f_{x}: \tilde{G}_{x} \xrightarrow{\cong}
\tilde{H}_{f(x)}$ induced by the local weak equivalence $f$ determine
a homotopy
\begin{equation*}
R(G) \times \underline{\mathbf{1}} \to \mathbf{Grp}(\mathcal{C})
\end{equation*}
from $F(G)$ to $F(H)\cdot f_{\ast}$. It follows that $F(G)$ and
$F(H)$ represent the same element of
$\pi_{0}H(\ast,\mathbf{Grp}(\mathcal{C}))$, and so the assignment $G
\mapsto [F(G)]$ induces a function
\begin{equation*}
\Phi: \pi_{0}\mathbf{Gerbe}(\mathcal{C}) \to \pi_{0}H(\ast,\mathbf{Grp}(\mathcal{C})).
\end{equation*}

Suppose that 
\begin{equation*}
\ast \xleftarrow{\simeq} A \xrightarrow{K} \mathbf{Grp}(\mathcal{C}) 
\end{equation*}
is a $2$-cocycle with coefficients in $\mathbf{Grp}(\mathcal{C})$.
Then $K$ consists of $2$-groupoid morphisms $K(U): A(U) \to
\mathbf{Grp}(\mathcal{C})(U)$, and hence induces composite morphisms
\begin{equation*}
A(U) \xrightarrow{K(U)} \mathbf{Grp}(\mathcal{C})(U)
\xrightarrow{ev_{U}} \mathbf{grp}.
\end{equation*}
Here, $ev_{U}: \mathbf{Grp}(\mathcal{C})(U) \to \mathbf{grp}$ is the
$2$-groupoid morphism
which is defined by $U$-sections.

Write $E_{A}K(U)$ for the Grothendieck construction corresponding to
the composite $ev_{U}K(U)$. Then the assignment $U \mapsto E_{A}K(U)$
defines a presheaf of groupoids $E_{A}K$. From 
Section 2, we see that there is a canonical morphism $p: E_{A}K \to
\pi A$ of presheaves of groupoids; it is defined in sections to be the
identity on objects, and it sends a class $[(f,\alpha)]$ to the class
$[\alpha]$.

\begin{lemma}\label{lem 16}
Suppose that $K: A \to \mathbf{Grp}(\mathcal{C})$ is a $2$-cocycle
over the terminal object $\ast$ taking values in sheaves of
groups. Then the presheaf of groupoids $E_{A}K$ is a gerbe.
\end{lemma}

\begin{proof}
The map
$\pi_{0}E_{A}K \to \pi_{0}(\pi A)$ is an isomorphism of presheaves, since
each $2$-functor $K(U)$ takes values
in groups.
\end{proof}

Take $i \in A(U)$ and let $K(i)$ be the corresponding sheaf of groups
on $\mathcal{C}/U$. The functor $\phi: K(i) \to p/i$ of presheaves
of groupoids on $\mathcal{C}/U$ is defined in sections corresponding
to an object $\psi: V \to U$ of $\mathcal{C}/U$ by sending the group
element $f \in K(i)(V)= K(\psi^{\ast}(i))(V)$ to the class
$[(f,1_{i})]$. 

\begin{lemma}\label{lem 17}
Suppose that $K: A \to \mathbf{Grp}(\mathcal{C})$ is a $2$-cocycle
over $\ast$ taking values in sheaves of groups, and choose $i \in
A(U)$. Then the homomorphism $\gamma_{i}: K(i) \to \hom(i,i) \subset
E_{A}K$ defined by $f \mapsto [(f,1_{i})]$ induces an isomorphism of
sheaves of groups on $\mathcal{C}/U$.
\end{lemma}

\begin{proof}
Suppose that $[(g,\alpha)]$ is an element of $\hom(i,i)$. By Lemma
\ref{lem 7} there is a covering sieve $R \subset \hom(\ ,U)$ such
that there is a $2$-cell $h_{\phi}: \phi^{\ast}(\alpha) \to
1_{\alpha^{\ast}(i)}$ for all $\phi \in R$. It follows that, locally,
$[(g,\alpha)]$ is in the image of $\gamma_{i}$.

Take group elements $f,g \in K(i)(U)$ and suppose that $\gamma_{i}(f)
= \gamma_{i}(g)$. Then there is a $2$-cell $h: 1_{i} \to 1_{i}$ in
$A(U)$ such that the diagram
\begin{equation*}
\xymatrix@R=6pt{
\ast \ar[dr]^{f} \ar[dd]_{h_{\ast}}  & \\
& \ast \\
\ast \ar[ur]_{g}
}
\end{equation*}
commutes in $K(i)(U)$. The presheaf of groupoids $A(i,i)$ is a
\v{C}ech object by Lemma \ref{lem 8} so that there is a covering
$\phi: V \to U$ such that $\phi^{\ast}(h) = 1$ for all members
$\phi$ of the cover. But then $\phi^{\ast}(h_{\ast}) = 1$ for all
$\phi$, and so $h_{\ast} = 1$ since $K(i)$ is a sheaf of groups.
\end{proof}

\begin{corollary}
Suppose that $K: A \to \mathbf{Grp}(\mathcal{C})$ is a $2$-cocycle
over $\ast$ taking values in sheaves of groups, and choose $i \in
A(U)$. Then the corresponding map $\gamma_{i}: K(i) \to \pi/i$ is a local
equivalence of presheaves of groupoids on $\mathcal{C}/U$.
\end{corollary}

\begin{proof}
The map $\gamma_{i}: K(i) \to \pi/i$ takes the group element $f$ to the
automorphism $[(f,1_{i})]$ of the object $[1_{i}]: \pi(i) \to i$. The
induced map $K(i) \to \hom([1_{i}],[1_{i}])$ is a surjection of
presheaves of groups. The composite
\begin{equation*}
K(i) \to \hom([1_{i}],[1_{i}]) \to \hom(i,i)
\end{equation*}
is locally monic by Lemma \ref{lem 17}. 

The map $\gamma_{i}: K(i) \to \pi/i$ is a sectionwise surjection on path
components.
\end{proof}

\begin{corollary}\label{cor 19}
Suppose that the diagram
\begin{equation*}
\xymatrix@R=6pt{
& A \ar[dd]_{\theta} \ar[dr]^-{K} \ar[dl]_{\simeq} & \\
\ast && \mathbf{Grp}(\mathcal{C}) \\
& B \ar[ul]^{\simeq} \ar[ur]_-{G}
}
\end{equation*}
is a morphisms of $2$-cocycles. Then the induced map $\theta: E_{A}K \to
E_{B}G$ is a local weak equivalence of presheaves of groupoids.
\end{corollary}

\begin{proof}
The diagram
\begin{equation*}
\xymatrix{
K(i) \ar[r]^-{\gamma_{i}} \ar[d]_{=} 
& \hom(i,i) \ar[d]^{\theta_{\ast}} \\
G(\theta(i)) \ar[r]_-{\gamma_{\theta(i)}} & \hom(\theta(i),\theta(i))
}
\end{equation*}
commutes, so that $\theta$ induces an isomorphism on all sheaves of
fundamental groups by Lemma \ref{lem 17}. The map $\theta$ induces an
isomorphism on sheaves of path components by Lemma \ref{lem 16}.
\end{proof}

It follows that the assignment $K \mapsto E_{A}K$ defines a functor
$H(\ast,\mathbf{Grp}(\mathcal{C})) \to \mathbf{Gerbe}(\mathcal{C})$ and hence a
function
\begin{equation*}
\Psi: \pi_{0}H(\ast,\mathbf{Grp}(\mathcal{C})) \to \pi_{0}(\mathbf{Gerbe}(\mathcal{C})).
\end{equation*}

\begin{theorem}\label{thm 20}
The functions $\Phi$ and $\Psi$ are inverse to each other, and define
a bijection
\begin{equation*}
\pi_{0}(\mathbf{Gerbe}(\mathcal{C})) \cong \pi_{0}H(\ast,\mathbf{Grp}(\mathcal{C})).
\end{equation*}
\end{theorem}

\begin{proof}
The relation $\Psi\Phi = 1$ is a consequence of Lemma \ref{lem 13}.

Suppose that $K: A \to \mathbf{Grp}(\mathcal{C})$ is a $2$-cocycle
over $\ast$. There is a $2$-groupoid morphism $\omega: A \to
R(E_{A}K)$ which is the identity on objects, sends the $1$-cell
$\alpha: i \to j$ to the $1$-cell $[(e,\alpha)]$, and sends the
$2$-cell $h: \alpha \to \beta$ to the $2$-cell
\begin{equation*}
[(h_{\ast},1)]: [(e,\alpha)] \to [(e,\beta)].
\end{equation*}

The presheaf of groupoids $\pi R(E_{A}K)$ has the same objects as $A$; it
is a \v{C}ech object (by Lemma \ref{lem 8}), in which there is a morphism $i \to j$ in
$R(E_{A}K)$ if and only if there is a $1$-cell $i \to j$ in
$A$. It follows that the morphism $\omega$ induces an isomorphism on
presheaves of path components.
It also follows
that the composite 
\begin{equation*}
A \xrightarrow{\omega} R(E_{A}K) \xrightarrow{F(E_{A}K)}
\mathbf{Grp}(\mathcal{C})
\end{equation*}
defines a group-valued $2$-cocycle on $\pi_{0}A$. This composite sends
the object $i \in A$ to the presheaf of groups $E_{A}K(i,i)$,
sends a $1$-cell $\alpha: i \to j$ to the homomorphism
$c_{\alpha}: E_{A}K(i,i) \to E_{A}K(j,j)$ which is defined by
conjugation with $[(e,\alpha)]$, and sends a $2$-cell $h: \alpha \to
\beta$ to the homotopy defined by conjugation with the element
$[(h_{\ast},1)]$.

The assignments $f \mapsto [(f,1)]$ define homomorphisms
\begin{equation*}
\gamma_{i}: K(i) \to E_{A}K(i,i).
\end{equation*}
which induce isomorphisms of associated sheaves, by Lemma \ref{lem 17}.
The morphisms $\gamma_{i}$ further determine a homotopy
\begin{equation*}
\gamma: A \times \underline{\mathbf{1}} \to \mathbf{Grp}(\mathcal{C}).
\end{equation*}
from the cocycle $K$ to the cocycle $F(E_{A}K)\omega$.
It follows that there is a path
\begin{equation*}
F(E_{A}K) \sim F(E_{A}K)\omega  \sim \gamma \sim K
\end{equation*}
in the cocycle category, and so $\Psi\Phi = 1$ as required. 
\end{proof}

\section{Homotopy classification of gerbes}

Suppose that $G$ is a presheaf of groupoids on $\mathcal{C}$, with
automorphism sheaves $\tilde{G}_{x}$, $x \in G(U)$. The presheaf of
$2$-groupoids $G_{\ast}$ has the same objects as $G$; the $1$-cells $x
\to y$ of $G_{\ast}(U)$ are the sheaf isomorphisms 
$\tilde{G}_{x} \to \tilde{G}_{y}$,
and the $2$-cells of $G_{\ast}(U)$ are the homotopies of
isomorphisms. There is a $2$-functor $\nu_{G}: G_{\ast} \to
\mathbf{Grp}(\mathcal{C})$ which is defined by sending $x$ to $\tilde{G}_{x}$,
and is the identity on sheaf isomorphisms and homotopies. The
canonical $2$-cocycle $F(G): R(G) \to \mathbf{Grp}(\mathcal{C})$
factors uniquely through a cocycle $F(G)_{\ast}: R(G) \to
G_{\ast}$ in the category of presheaves of $2$-groupoids.

Suppose that
$\mathcal{F} \subset \mathbf{Grp}(\mathcal{C})$ is a subobject of
$\mathbf{Grp}(\mathcal{C})$ such that
\begin{itemize}
\item[1)] 
the imbedding is full: all simplicial presheaf maps
\begin{equation*}
\mathcal{F}(H,K) \to \mathbf{Grp}(\mathcal{C})(H,K) 
\end{equation*}
are isomorphisms,
\item[2)] 
$\mathcal{F}$ is a presheaf of groupoids, so that all
classes $\Ob(\mathcal{F})(U)$ are sets,
\end{itemize} 
We shall say that a subobject $\mathcal{F}$ of the diagram of
$2$-groupoids $\mathbf{Grp}(\mathcal{C})$ which satisfies these
conditions is a {\it full subpresheaf} of $\mathbf{Grp}(\mathcal{C})$.

The image of the presheaf of $2$-groupoids $G_{\ast}$ in
$\mathbf{Grp}(\mathcal{C})$ which arises from a presheaf of groupoids
$G$ is an example of such an object $\mathcal{F}$.

\begin{lemma}\label{lem 21}
Suppose that $\mathcal{F} \subset \mathcal{F}'$ are full subpresheaves
of $\mathbf{Grp}(\mathcal{C})$. Suppose further that every
automorphism group $\mathcal{F}'_{x}$ of $\mathcal{F}'$ is locally
isomorphic to automorphism groups of $\mathcal{F}$. Then the inclusion
$\mathcal{F} \subset \mathcal{F}'$ is a local weak equivalence of
presheaves of $2$-groupoids.
\end{lemma}

\begin{proof}
Write $\alpha: \mathcal{F} \subset \mathcal{F}'$ for the inclusion
morphism.  Then $\alpha$ is full, and therefore induces a presheaf
monomorphism $\pi_{0}G_{\ast} \to \pi_{0}\mathcal{F}$. Every sheaf of
groups $\mathcal{F}'_{x} \in \mathcal{F}'(U)$ is locally isomorphic to
objects in the image of $\alpha$, by definition, so that
$\pi_{0}\mathcal{F} \to \pi_{0}\mathcal{F}'$ is a local epimorphism.

The assertion that $\alpha$ induces an isomorphism in all possible
sheaves of higher homotopy groups is a consequence of the fullness and
Lemma \ref{lem 1}.
\end{proof}

Say that two gerbes $G$ and $H$ are {\it locally equivalent} if there
is a covering family $U \to \ast$, $U \in \Ob(\mathcal{C})$, such that
the restricted gerbes $G\vert_{U}$ and $H\vert_{U}$ are locally weakly
equivalent on $\mathcal{C}/U$ for each object $U$ in the covering of
the terminal object $\ast$. If there is a local weak equivalence $G
\to H$ then $G$ and $H$ are locally equivalent in the sense just
described, but the converse is not true.

\begin{example}
Suppose that the presheaf of $2$-groupoids $\mathcal{F}$ is a full subpresheaf of $\mathbf{Grp}(\mathcal{C})$, and that there is a
$2$-cocycle 
\begin{equation*}
\ast \xleftarrow{\simeq} A \xrightarrow{F} \mathcal{F} \subset
\mathbf{Grp}(\mathcal{C})
\end{equation*}
over the terminal sheaf $\ast$. 

There is a covering family $U \to \ast$, $U \in \mathcal{C}$, such
that $A(U) \ne \emptyset$. In effect, $\Ob(A) \to \ast$ is a local
epimorphism, so there is a covering $U \to \ast$ such that there are
liftings
\begin{equation*}
\xymatrix{
& \Ob(A) \ar[d] \\
U \ar[r] \ar[ur]^{x_{U}} & \ast
}
\end{equation*}
where $x_{U}$ represents an object of $A(U)$. The presheaf of
groupoids $E_{A}F$ is locally connected by Lemma \ref{lem 16}, and the
maps 
\begin{equation*}
F(x_{U}) \to \hom_{E_{A}F}(x_{U},x_{U})
\end{equation*}
induce isomorphisms of associated sheaves of groups on $\mathcal{C}/U$
by Lemma \ref{lem 17}. It follows that the automorphism groups of the
Grothendieck construction $E_{A}F$ are locally equivalent to objects
of $\mathcal{F}$.
\end{example}

Write $\mathcal{F}-\mathbf{Gerbe}$ for the full subcategory of the
category of gerbes whose automorphism groups are locally equivalent to
sheaves of groups in $\mathcal{F}$.  The assignment $F \mapsto E_{A}F$
for a cocycle $F: A \to \mathcal{F}$ takes values in
$\mathcal{F}$-gerbes, so that there is a commutative diagram
\begin{equation*}
\xymatrix{
\pi_{0}H(\ast,\mathcal{F}) \ar[r] \ar[d] 
& \pi_{0}H(\ast,\mathbf{Grp}(\mathcal{C})) \ar[d]^{\cong} \\
\pi_{0}(\mathcal{F}-\mathbf{Gerbe}) \ar[r] & \pi_{0}(\mathbf{Gerbe})
}
\end{equation*}
Note that if $f: G \to H$ is a local weak equivalence of gerbes, then
$G$ is an $\mathcal{F}$-gerbe if and only if $H$ is an
$\mathcal{F}$-gerbe, and it follows that the induced map
\begin{equation*}
\pi_{0}(\mathcal{F}-\mathbf{Gerbe}) \to \pi_{0}(\mathbf{Gerbe})
\end{equation*}
is an injection.

\begin{theorem}\label{thm 23}
Suppose that $\mathcal{F}$ is a full subpresheaf of the
$\mathbf{Grp}(\mathcal{C})$. Then the Grothendieck construction defines a 
function
\begin{equation*}
\pi_{0}H(\ast,\mathcal{F}) \to \pi_{0}(\mathcal{F}-\mathbf{Gerbe}) 
\end{equation*}
which is a bijection.
\end{theorem}

\begin{proof}
Suppose given cocycles $F: A \to \mathcal{F}$ and $G: B \to \mathcal{F}$
such that $F$ and $G$ are in the same path component as cocycles
taking values in $\mathbf{Grp}(\mathcal{C})$. Then there is a string
of maps of cocycles
\begin{equation}\label{eq 2}
F = F_{0} \leftrightarrow F_{1} \leftrightarrow \dots \leftrightarrow F_{n} = G
\end{equation}
where $F_{i}: A_{i} \to \mathbf{Grp}(\mathcal{C})$ are cocycles in
$\mathbf{Grp}(\mathcal{C})$.

Suppose that $F: A \to \mathbf{Grp}(\mathcal{C})$ is a cocycle taking
values in sheaves of groups locally isomorphic to objects of $\mathcal{F}$
and that
\begin{equation*}
\xymatrix@C=6pt{
A \ar[rr]^{\alpha} \ar[dr]_{F} && A' \ar[dl]^{F'} \\
& \mathbf{Grp}(\mathcal{C})
}
\end{equation*}
is a morphism of $H(\ast,\mathbf{Grp}(\mathcal{C}))$. Take $x \in
A'(U)$. Then there is a covering family $\phi: V \to U$ with $1$-cells
$\phi^{\ast}(x) \to \alpha(y_{V})$ in $A'(V)$ for all $\phi$. It
follows that the group $F'(x)$ is locally isomorphic to groups of the
form $F(y_{V})$, and all of these are locally isomorphic to
objects of $\mathcal{F}$. Thus, the cocycle $F'$ takes values in
sheaves of groups locally isomorphic to objects of $\mathcal{F}$

It follows that all cocycles $F_{i}$ in the list (\ref{eq 2}) take
values in groups locally isomorphic to objects of $\mathcal{F}$.
Write $\mathcal{F}'$ for the presheaf of $2$-groupoids which is the
full subobject of $\mathbf{Grp}(\mathcal{C})$ on the sheaves of groups
appearing in the sets $\mathcal{F}(U)$ and all $F_{i}(\Ob(A_{i}))(U)$.
Then $\mathcal{F} \subset \mathcal{F}'$, and
Lemma \ref{lem 21} implies that this map of presheaves of
$2$-groupoids is a weak equivalence. The string of cocycles $F_{i}$ in
(\ref{eq 2}) all take values in $\mathcal{F}'$ by construction, and the map
\begin{equation*}
\pi_{0}H(\ast,\mathcal{F}) \to \pi_{0}H(\ast,\mathcal{F}')
\end{equation*}
is a bijection. It follows that the original cocycles $F$ and $G$ are
in the same path component of $H(\ast,\mathcal{F})$. The function
\begin{equation*}
\pi_{0}H(\ast,\mathcal{F}) \to \pi_{0}H(\ast,\mathbf{Grp}(\mathcal{C}))
\end{equation*}
is therefore a monomorphism, as is the function
\begin{equation*}
\pi_{0}H(\ast,\mathcal{F}) \to \pi_{0}(\mathcal{F}-\mathbf{Gerbe}).
\end{equation*}

Suppose that the every automorphism presheaf of the gerbe $H$ is
locally equivalent to an object of $\mathcal{F}$. Then all
automorphism sheaves of $H$ are locally isomorphic to automorphism
sheaves of $G$. Choose a full subpresheaf $\mathcal{F}' \subset
\mathbf{Grp}(\mathcal{C})$ whose $0$-cells are sheaves of groups
locally equivalent to objects of $\mathcal{F}$ and which contains both
$\mathcal{F}$ and $H_{\ast}$. Then the canonical cocycle
\begin{equation*}
F(H): R(H) \to \mathbf{Grp}(\mathcal{C}) 
\end{equation*}
takes values in $\mathcal{F}'$. The map $\mathcal{F} \to \mathcal{F}'$
is a local weak equivalence, and so $F(H)$ can be represented by a
cocycle taking values in $\mathcal{F}$.  It follows that the function
$\pi_{0}H(\ast,\mathcal{F}) \to \pi_{0}(\mathcal{F}-\mathbf{Gerbe})$
is surjective.
\end{proof}

\begin{corollary}\label{cor 24}
Suppose that $\mathcal{F}$ is a full subpresheaf of $2$-groupoids in $\mathbf{Grp}(\mathcal{C})$. Then there is a bijection
\begin{equation*}
[\ast,dB\mathcal{F}] \cong \pi_{0}(\mathcal{F}-\mathbf{Gerbe}).
\end{equation*}
\end{corollary}

Suppose that $G$ is a gerbe, and write $G-\mathbf{Gerbe}$ for the
category of gerbes which are locally equivalent to $G$. This category
coincides with the category $G_{\ast}-\mathbf{Gerbe}$ arising from the
full subpresheaf of $2$-groupoids $G_{\ast}$, and so we have the
following:

\begin{corollary}\label{cor 25}
Suppose that $G$ is a gerbe on a site $\mathcal{C}$ with associated
$2$-groupoid object $G_{\ast}$ of isomorphisms and homotopies of
automorphism sheaves of $G$. Then there is a bijection
\begin{equation*}
[\ast,dBG_{\ast}] \cong \pi_{0}(G-\mathbf{Gerbe}).
\end{equation*}
\end{corollary}

A special case of Corollary \ref{cor 25}, corresponding to the case of
a sheaf of groups $G$, was proved by Breen in \cite{Br}.

\begin{remark}
Recall that gerbes on $\mathcal{C}$ can be identified with gerbes
on the fibred site $\mathcal{C}/E$ up to natural equivalence. Given an
gerbe $G$, write $G_{E}$ for the corresponding gerbe on
$\mathcal{C}/E$. Then Corollary \ref{cor 25} gives a homotopy
classification
\begin{equation*} 
[\ast,dBG_{E\ast}] \cong \pi_{0}(G_{E}-\mathbf{Gerbe}).
\end{equation*}
for gerbes, up to local equivalence defined on the site $\mathcal{C}/E$.
\end{remark}

Suppose that $\mathcal{F}'$ is a full subpresheaf of
$\mathbf{Grp}(\mathcal{C})$, and write $\mathbf{Gerbe}(\mathcal{F'})$
for the full subcategory of gerbes $G$ such that $G_{\ast} \subset
\mathcal{F}'$ --- say that the objects of
$\mathbf{Gerbe}(\mathcal{F}')$ are the gerbes in $\mathcal{F'}$. The
category $\mathcal{F}-\mathbf{Gerbe}$ of gerbes with automorphism
sheaves locally isomorphic to objects of $\mathcal{F}$ is a filtered
colimit of subcategories $\mathbf{Gerbe}(\mathcal{F}')$, indexed over
all inclusions $\mathcal{F} \subset \mathcal{F}'$ of full
subpresheaves of $\mathbf{Grp}(\mathcal{C})$ such that every object of
$\mathcal{F'}$ is locally isomorphic to objects of $\mathcal{F}$. It
follows that there is an isomorphism
\begin{equation*}
\pi_{0}(\mathcal{F}-\mathbf{Gerbe}) \cong \varinjlim_{\mathcal{F}
\underset{\simeq}{\subset} \mathcal{F}'} \pi_{0}(\mathbf{Gerbe}(\mathcal{F}')).
\end{equation*}

Write $\St(\mathcal{F})$ and $\St(\pi\mathcal{F})$ for the stack
completions (really, fibrant models) for the presheaf of $2$-groupoids
$\mathcal{F}$ and its path component object $\pi \mathcal{F}$. The
path component object is a groupoid of outer automorphisms and its
stack completion $\St(\pi\mathcal{F})$ is the stack of {\it bands}
({\it liens}) for $\mathcal{F}$. These stack completion constructions
are functorial, since the underlying model structures are cofibrantly
generated \cite{L}.

A (global) {\it band} $L$ is a global section of the presheaf of
groupoids $\St(\pi\mathcal{F})$, or equivalently \cite{J3} a torsor
for the presheaf of outer automorphism groupoids
$\pi\mathcal{F}$. 

Write $p_{\mathcal{F}}$ for the composite
\begin{equation*}
\mathcal{F} \to \pi\mathcal{F} \to \St(\pi\mathcal{F}).
\end{equation*}
The homotopy fibre over a global band $L$ of the induced map $B\mathcal{F} \to
B\St(\pi\mathcal{F})$ is the classifying object $B(p_{\mathcal{F}}/L)$ of the
simplicial groupoid $p_{\mathcal{F}}/L$ \cite{J2}.

The objects of $2$-cocycle category
$H(\ast,B(p_{\mathcal{F}}/L))$ can be identified with the collection of pairs
$(\nu,\phi)$ consisting of 
a $2$-cocycle
\begin{equation*}
\ast \xleftarrow{\simeq} A \xrightarrow{\phi} \mathcal{F}
\end{equation*}
and a natural (iso)morphism $\nu: \phi_{\ast} \to L$ in
$\St(\pi\mathcal{F})$, where $\phi_{\ast}: \pi(A) \to
\St(\pi\mathcal{F})$ is the unique induced morphism in the diagram
\begin{equation}\label{eq 3}
\xymatrix{
A \ar[r]^{\phi} \ar[d] & \mathcal{F} \ar[d] \\
\pi(A) \ar[r]_{\phi_{\ast}} & \St(\pi\mathcal{F})
}
\end{equation}
and $L$ has been identified with the composite
\begin{equation*}
\pi(A) \to \ast \xrightarrow{L} \St(\pi\mathcal{F}).
\end{equation*}
 The morphisms $f: (\phi,\nu) \to (\phi',\nu')$ are cocycle morphisms
\begin{equation*}
\xymatrix@R=8pt{
A \ar[dd]_{f} \ar[dr]^{\phi} & \\
& \mathcal{F} \\
A' \ar[ur]_{\phi'}
}
\end{equation*}
such that $\nu' \cdot f_{\ast} = \nu: \phi_{\ast} \to L$. 

Suppose that $G$ is a gerbe in $\mathcal{F}$.
Consider the diagram
\begin{equation}\label{eq 4}
\xymatrix{
G \ar[r] \ar[dr] & R(G) \ar[r]^{F(G)} \ar[d] & \mathcal{F} \ar[r] \ar[d]
& \St(\mathcal{F}) \ar[d] \\
& \pi R(G) \ar[r]_-{F(G)_{\ast}} & \pi\mathcal{F} \ar[r]_-{\eta} 
& \St(\pi\mathcal{F})
}
\end{equation}
and write $\omega_{G} = \eta F(G)_{\ast}: \pi R(G) \to \St(\pi\mathcal{F})$.

An $L$-gerbe $(G,\nu_{G})$ in $\mathcal{F}$ is a gerbe $G$ in
$\mathcal{F}$ together with a natural isomorphism $\nu_{G}:
\omega_{G} \to L$, where $L$ has been identified with the composite
\begin{equation*}
G \to \ast \xrightarrow{L} \St(\pi\mathcal{F}).
\end{equation*}
Observe that there is a canonical natural isomorphism
\begin{equation*}
h_{f\ast}: F(G)_{\ast} \xrightarrow{\cong} F(H)_{\ast}f_{\ast}
\end{equation*}
for any morphism $f: G \to H$ of gerbes.  A morphism $f: (G,\nu_{G})
\to (H,\nu_{H})$ of $L$-gerbes is a morphism $f: G \to H$ of gerbes
such that the diagram of natural isomorphisms
\begin{equation}\label{eq 5}
\xymatrix{
\omega_{G} \ar[d]_{\nu_{G}} \ar[r]^-{\eta(h_{f\ast})}
& \omega_{H}f_{\ast} \ar[d]^{\nu_{H}f_{\ast}} \\
L \ar[r]_{1} & L
}
\end{equation}
commutes. The natural isomomorphisms $h_{f\ast}$ arising from gerbe
morphisms $f$ are coherent; this gives the law of composition for a
category of $L$-gerbes in $\mathcal{F}$, which will be denoted by
$\mathbf{Gerbe}(\mathcal{F})/L$.

Every $L$-gerbe $(G,\nu_{G})$ in $\mathcal{F}$ determines an object
$(F(G),\nu_{G})$ in the cocycle category $H(\ast,B(p_{\mathcal{F}}/L))$.

Suppose that $f: (G,\nu_{G}) \to (H,\nu_{H})$ is a morphism of
$L$-gerbes in $\mathcal{F}$. Then the homotopy of cocycles $h_{f}: R(G) \times
\underline{\mathbf{1}} \to \mathcal{F}$ from $F(G)$ to $F(H)f_{\ast}$
determines a diagram of path component groupoid morphisms
\begin{equation*}
\xymatrix{
\pi R(G) \ar[d] \ar[drr]^{F(G)_{\ast}} &&& \\
\pi R(G) \times \underline{\mathbf{1}} \ar[rr]^{h_{f\ast}} && \pi \mathcal{F} \ar[r]^-{\eta}  
& \St(\pi\mathcal{F}) \\
\pi R(G) \ar[u] \ar[r]_{f_{\ast}} & \pi R(H) \ar[ur]_{F(H)_{\ast}}
}
\end{equation*}
The natural isomorphism $\nu_{G}: \eta F(G)_{\ast} \to L$ extends
uniquely to a natural isomorphism $\nu_{h}: \eta h_{f\ast} \to L$,
and $\nu_{h}$ restricts to $\nu_{H}f_{\ast}: \eta F(H)_{\ast} \to L$
on $\pi R(G) \times \{ 1 \}$ on account of the commutativity of the
diagram (\ref{eq 5}). 

It follows that every morphism $f: (G,\nu_{G}) \to (H,\nu_{H})$ of
$L$-gerbes determines a path between the associated objects
$(F(G),\nu_{G})$, $(F(H),\nu_{H})$ in the cocycle category, and that
there is a function
\begin{equation*}
\Phi_{\mathcal{F}}: \pi_{0}(\mathbf{Gerbe}(\mathcal{F})/L) 
\to \pi_{0}H(\ast,B(p_{\mathcal{F}}/L))
\end{equation*}
which is defined by $\Phi([(G,\nu_{G})]) = [(F(G),\nu_{G})]$.

Suppose that the $2$-cocycle
\begin{equation*}
\xymatrix{
\ast & A \ar[l]_{\simeq} \ar[r]^{\phi} \ar[d] & \mathcal{F} \ar[d] \\
& \pi(A) \ar[r]_{\phi_{\ast}} & \St(\pi\mathcal{F})
}
\end{equation*}
and the natural isomorphism $\nu: \phi_{\ast} \to L$ in
$\St(\pi\mathcal{F})$ define an object $(\nu,\phi)$ of the cocycle
category $H(\ast,B(p_{\mathcal{F}}/L))$. Then the associated presheaf
of groupoids $E_{A}\phi$ is a gerbe which has automorphism sheaves
locally isomorphic to objects of $\mathcal{F}$, and then from
the proof of Theorem \ref{thm 20} we know that there is a
homotopy
\begin{equation*}
\xymatrix{
A \ar[d] \ar[drr]^{\phi} && \\
A \times \underline{\mathbf{1}} \ar[rr]^{\gamma} 
&& \mathcal{F}' \\
A \ar[u] \ar[r]_-{\omega} & R(E_{A}\phi) \ar[ur]_{F(E_{A}\phi)}
}
\end{equation*}
where $\mathcal{F}'$ is a full subpresheaf of
$\mathbf{Grp}(\mathcal{C})$ containing $\mathcal{F}$ such that the map
$\mathcal{F} \subset \mathcal{F}'$ is a local weak equivalence.  We
also know that the induced map $\omega_{\ast}: \pi A \to \pi
R(E_{A}\phi)$ is an isomorphism. It follows that the induced natural
isomorphism (or homotopy)
\begin{equation*}
\gamma_{\ast}: \phi_{\ast} \xrightarrow{\cong} F(E_{A}\phi)_{\ast}\omega_{\ast}
\end{equation*}
of functors $\pi A \to \St(\mathcal{F}')$ induces a unique natural isomorphism
\begin{equation*}
\eta F(E_{A}\phi)_{\ast} \xrightarrow{\tilde{\nu}} L
\end{equation*}
which restricts to the isomorphism $\nu: \phi_{\ast} \to L$ along the
homotopy 
\begin{equation*}
\pi (A) \times \underline{\mathbf{1}} \xrightarrow{\gamma_{\ast}}
\pi \mathcal{F}' \to \St(\pi\mathcal{F}').
\end{equation*}
In other words, $(E_{A}\phi,\tilde{\nu})$ is an $L$-gerbe in $\mathcal{F}'$.

Suppose that $f: (\phi,\nu) \to (\phi',\nu')$ is a morphism of the
$2$-cocycle category $H(\ast,B(p_{\mathcal{F}}/L))$. Then there is a full
subpresheaf $\mathcal{F}'' \subset \mathbf{Grp}(\mathcal{C})$
containing $\mathcal{F}$, such that $\mathcal{F} \subset
\mathcal{F}''$ is a weak equivalence and such that the associated
gerbes $E_{A}\phi$ and $E_{A'}\phi'$ are gerbes in
$\mathcal{F}''$. There is a diagram of homotopies
\begin{equation*}
\xymatrix{
\phi_{\ast}(i) \ar[r]^-{=} \ar[d]_{\gamma_{i}} 
& \phi'_{\ast}(f(i)) \ar[d]^{\gamma_{f(i)}} \\
\eta \Aut(i) \ar[r]_{f_{\ast}} & \eta \Aut(f(i))
}
\end{equation*}
where $\Aut(i)$ is the sheaf of automorphisms of $i$ in $E_{A}\phi$
and $\Aut(f(i))$ is the sheaf of automorphisms of $f(i)$ in
$E_{A'}\phi'$. Then the morphisms $\nu_{i}: \phi(i) \to L$ and
$\nu'_{f(i)} : \phi'(f(i)) \to L$ coincide on $\phi(i) = \phi'(f(i))$
since $f$ is a morphism of the cocycle category
$H(\ast,B(p_{\mathcal{F}}/L))$. Furthermore, the vertical isomorphisms uniquely
determine the natural isomorphisms $\tilde{\nu}: \eta \Aut(i) \to L$
and $\tilde{\nu}'f_{\ast}: \eta \Aut(f(i)) \to L$, respectively. It
follows that the map $f_{\ast}: E_{A}\phi \to E_{A'}\phi'$ defines a
morphism of $L$-gerbes in $\mathcal{F}''$. We therefore have a well
defined function
\begin{equation*}
\Psi: \pi_{0}H(\ast,B(p_{\mathcal{F}}/L)) \to \varinjlim_{\mathcal{F} 
\underset{\simeq}\subset \mathcal{F}'} 
\pi_{0}(\mathbf{Gerbe}(\mathcal{F}')/L).
\end{equation*}

\begin{theorem}\label{thm 27}
Suppose that $L$ is a band.
Then the function
\begin{equation*}
\end{equation*}
is a bijection.
\end{theorem}

\begin{proof}
Suppose that $\mathcal{F} \subset \mathcal{F'}'$ is a weak equivalence
of full subpresheaves of $\mathbf{Grp}(\mathcal{C})$. Then the diagram
\begin{equation*}
\xymatrix{
\pi_{0}(\mathbf{Gerbe}(\mathcal{F})/L) \ar[r]^{\Phi_{\mathcal{F}}} \ar[d] 
& \pi_{0}H(\ast,B(p_{\mathcal{F}}/L)) \ar[d]^{\cong} \\
\pi_{0}(\mathbf{Gerbe}(\mathcal{F'})/L) \ar[r]_{\Phi_{\mathcal{F}'}} 
& \pi_{0}H(\ast,B(p_{\mathcal{F}'}/L))
}
\end{equation*}
commutes, where the indicated vertical map is a bijection since the
comparison map $B(p_{\mathcal{F}}/L)) \to B(p_{\mathcal{F}'}/L)$ is a
local weak equivalence. It follows that the maps $\Phi_{\mathcal{F}'}$
together induce a function
\begin{equation*}
\Phi: \varinjlim_{\mathcal{F} \underset{\simeq}{\subset} \mathcal{F}'} 
\pi_{0}(\mathbf{Gerbe}(\mathcal{F}')/L) \to 
\pi_{0}H(\ast,B(p_{\mathcal{F}}/L)).
\end{equation*}
The function $\Phi$ is the inverse of $\Psi$.
\end{proof}

\begin{corollary}\label{cor 28}
Suppose that $L$ is a band. Then there are bijections
\begin{equation*}
[\ast,B(p_{\mathcal{F}}/L)] \cong \pi_{0}H(\ast,B(p_{\mathcal{F}}/L))
\cong \varinjlim_{\mathcal{F} \underset{\simeq}\subset \mathcal{F}'}
\pi_{0}(\mathbf{Gerbe}(\mathcal{F}')/L).
\end{equation*}
\end{corollary}
\vfill\eject


\catcode`\@=11 \renewenvironment{thebibliography}[1]{
\@xp\section\@xp*\@xp{\refname}%
\normalfont\footnotesize\labelsep
.5em\relax\renewcommand\theenumiv{\arabic{enumiv}}\let\p@enumiv\@empty
\list{\@biblabel{\theenumiv}}{\settowidth\labelwidth{\@biblabel{#1}}%
\leftmargin\labelwidth \advance\leftmargin\labelsep
\usecounter{enumiv}}%
\sloppy \clubpenalty\@M
\widowpenalty\clubpenalty \sfcode`\.=\@m }

\def\@biblabel#1{\@ifnotempty{#1}{[#1]}}
\catcode`\@=\active


\bibliographystyle{amsalpha}

\enddocument